\documentclass[11pt, eqno]{article}

\usepackage{amsmath,amssymb,latexsym}

\newcommand{\N}{{\mathbb N}}
\newcommand{\C}{{\mathbb C}}
\newcommand{\R}{{\mathbb R}}

\newcommand{\Z}{{\mathbb Z}}
\newcommand{\CP}{{\mathbb CP}}
\newcommand{\RP}{{\mathbb RP}}

\newtheorem{theorem}{Theorem}[section]

\newtheorem{remark}[theorem]{Remark}

\newtheorem{lemma}[theorem]{Lemma}

\newtheorem{proposition}[theorem]{Proposition}

\newfont{\Bb}{msbm10 scaled\magstephalf}

\begin{document}
\title{Symplectic fixed points and Lagrangian intersections on \\ weighted projective spaces}
\author{Guangcun Lu
\thanks{Partially supported by the NNSF 19971045  and 10371007 of China.}\\
Department of Mathematics,  Beijing Normal University\\
Beijing 100875,   P.R.China\\
(gclu@bnu.edu.cn)}

\date{January 10, 2006}
\maketitle \vspace{-0.1in}

\abstract{In this note we observe that Arnold conjecture for the
Hamiltonian maps still holds on weighted projective spaces
$\CP^n({\bf q})$, and that Arnold conjecture for the Lagrange
intersections for $(\CP^n({\bf q}), \RP^n({\bf q}))$ is also true
if each weight $q_i\in {\bf q}=\{q_1,\cdots, q_{n+1}\}$ is odd.}
 \vspace{-0.1in}
\medskip

\section{Introduction}

A famous conjecture by Arnold \cite{Ar} claimed that every exact
symplectic diffeomorphism on a closed symplectic manifold
$(P,\omega)$ has at least as many fixed points as the critical
points of a smooth function on $P$. The homological form of it can
be stated as: For a Hamiltonian map $\phi$, i.e., a time $1$-map
of a time-dependent Hamiltonian vector field $X_{h_t}$, $0\le t\le
1$,  the number of fixed points of $\phi$ satisfies the estimates
$$
(AC_1) \left\{\!\!\!\!\!\!\begin{array}{ll}
&\sharp{\rm Fix}(\phi)\ge CL(P)+1,\\
&\sharp{\rm Fix}(\phi)\ge SB(P)\quad{\rm if}\;{\rm every}\;{\rm
point}\;{\rm of}\;{\rm Fix}(\phi)\;{\rm is}\;{\rm nondegenerate}.
\end{array}
\right.
$$
Here $CL(P)$ is the cuplength of $P$ and $SB(P)$ is the sum of the
Betti numbers of $P$. More generally, for a closed Lagrange
submanifold $L$ in $(P,\omega)$ Arnold also conjectured:
$$
(AC_2) \left\{\!\!\!\!\!\begin{array}{ll}
&\sharp(L\cap(\phi(L))\ge CL(L)+1,\\
&\sharp(L\cap(\phi(L))\ge SB(L)\quad{\rm if}\;L\pitchfork\phi(L).
\end{array}
\right.
$$
After Conley and Zehnder \cite{CoZe} first proved $(AC_1)$ for the
standard symplectic torus $T^{2n}$, Fortune showed that $(AC_1)$
holds on $\CP^n$ with the standard structure. By generalizing the
idea of Gromov \cite{Gr}, Floer \cite{Fl1}-\cite{Fl3} found a
powerful approach to prove $(AC_1)$ and $(AC_2)$ for a large class
of symplectic manifolds and their Lagrangian submanifolds; also
see \cite{Ho} for a different method in the case $\pi_2(P, L)=0$.
$(AC_2)$ was proved for $(\CP^n,\RP^n)$ in \cite{ChJi} and
\cite{Gi}. Furthermore generalizations for $(AC_2)$ was made by Oh
in \cite{Oh1}-\cite{Oh3}.  Recently, by furthermore developing
Floer's method Fukaya-Ono \cite{FuO} and Liu-Tian \cite{LiuT}
proved the second claim in $(AC_1)$ for all closed symplectic
manifolds. The obstruction theory for Lagrangian intersection was
developed in \cite{FuOOO} very recently, and more general results
for $(AC_2)$ was also obtained.  For more complete history and
references of the conjectures we refer to  \cite{FuOOO},
\cite{McSa} and \cite{Se}.

 A symplectic orbifold is a natural generalization of a symplectic manifold.
Recall that a symplectic orbifold is a pair $(M,\omega)$
consisting of an orbifold $M$ and a closed non-degenerate $2$-form
$\omega$ on it. That is, $\omega$ is a differential form which in
each local representation is a closed nondegenerate $2$-form. Many
definitions on symplectic manifolds, e.g., Hamiltonian maps,
symplectic group actions, moment maps and Hamiltonian actions can
carry over verbatim to the category of symplectic orbifolds, cf.,
\cite{LeTo}. Ones can, of course, raise the corresponding ones of
the Arnold conjectures above on closed symplectic orbifolds.

Weighted projective spaces are typical symplectic orbifolds. Let
${\bf q}=(q_1,\cdots,q_{n+1})$ be a $(n+1)$-tuple of positive
integers. Recall that the weighted (twisted) projective space of
type ${\bf q}$ is defined by
$$
\CP^n({\bf
q})=(\C^{n+1}\setminus\{0\})/\C^\ast,
$$
where $\C^\ast=\C\setminus\{0\}$ acts on $\C^{n+1}\setminus\{0\}$
by
\begin{equation}\label{e:1.1} \alpha\cdot {\bf
z}=(\alpha^{q_1}z_1,\cdots, \alpha^{q_{n+1}}z_{n+1})
\end{equation}
 for ${\bf z}=(z_1,\cdots, z_{n+1})\in\C^{n+1}$ and $\alpha\in\C^\ast$.
 Note that the above $\C^\ast$-action is free iff $q_i=1$ for every
$i=1,\cdots,n+1$. If the largest common divisor $lcd(q_1,\cdots,
q_{n+1})=1$, $\CP^n({\bf q})$ has only isolated orbifold
singularities.
 Let $[{\bf z}]_{\bf q}$ denote the orbit of
${\bf z}\in\C^{n+1}\setminus\{0\}$ under the above
$\C^\ast$-action, i.e.,  a point in $\CP^n({\bf q})$. Denote by
$m({\bf z})$ the largest common divisor of the set $\{q_j\,|\,
z_j\ne 0\}$. The orbifold structure group $\Gamma_{[{\bf z}]_{\bf
q}}$ of $[{\bf z}]_{\bf q}$ is isomorphic to $\Z/m\Z$. So $[{\bf
z}]_{\bf q}$ is a smooth point of $\CP^n({\bf q})$ if and only if
$m({\bf z})=1$. Clearly,
 each point $[{\bf z}]_{\bf q}\in\CP^n({\bf q})$ with all $z_j\ne 0$, is a smooth
point. As on usual complex projective spaces ones can use
symplectic reduction to describe the symplectic orbifold structure
on $\CP^n({\bf q})$. Indeed, as showed in Proposition 2.8 of [Go]
the action of $S^1\equiv\R/2\pi\Z$,
\begin{equation}\label{e:1.2}
A^{\bf q}_s(z_1,\cdots, z_{n+1})=(e^{iq_1s}z_1,\cdots,
e^{iq_{n+1}s}z_{n+1})\;\forall s\in\R,
\end{equation}
is a Hamiltonian circle action on $(\C^{n+1},\omega_0)$ with a
moment map
\begin{equation}\label{e:1.3}
K_{\bf q}(z_1,\cdots,
z_{n+1})=\frac{1}{2}\sum^{n+1}_{i=1}q_i|z_i|^2;
\end{equation}
and each $t\ne 0$ is a regular value of $K_{\bf q}$. The circle
action on $K_{\bf q}^{-1}(t)$ is locally free, and thus
$\CP^n({\bf q})\cong K_{\bf q}^{-1}(t)/S^1$ for each $t\ne 0$.
Denote by $S^{2n+1}({\bf q})=K_{\bf q}^{-1}(\frac{1}{2})$, and by
$\Pi:S^{2n+1}({\bf q})\to S^{2n+1}({\bf q})/S^1=\CP^n({\bf q})$ be
the natural projection.  The reduction symplectic form
$\omega_{\rm FS}^{\bf q}$ on $S^{2n+1}({\bf q})/S^1$ is called
{\it standard orbifold symplectic form} on $\CP^n({\bf q})$
without special statements (since $\omega_{\rm FS}^{\bf
q}=\omega_{\rm FS}$ has integral one on $\CP^1\subset\CP^n$ if
each component $q_i$ is equal to $1$ in ${\bf q}$). It was
computed in \cite{Ka} that $H^i(\CP^n({\bf q});\Z)=\Z$
 for $i=2k$ and $0\le k\le n$, and zero for other $i$.
Let  $\gamma_k$ denote the canonical generator of the group
$H^{2k}(\CP^n({\bf q});\Z)$ for $1\le k\le n$. It was also proved
in \cite{Ka} that the multiplication is given by
$$
\gamma_k\gamma_j=\frac{l_k^{\bf q}l_j^{\bf q}}{l_{k+j}^{\bf
q}}\gamma_{k+j}\;{\rm if}\;k+j\le n,\quad{\rm and}\;
\gamma_k\gamma_j=0\quad {\rm if}\; k+j\ge n.
$$
Here for $1\le k\le n$,
$$
l^{\bf q}_k={\rm lcm}\biggl(\biggl\{\frac{q_{i_1}\cdots
q_{i_{k+1}}}{{\rm gcd}(q_{i_1},\cdots, q_{i_{k+1}})}\biggm| 1\le
i_1<\cdots<i_{k+1}\le n+1\biggr\}\biggr),
$$
and $l_k^{\bf q}l_j^{\bf q}/l_{k+j}^{\bf q}$ is always an integer.
It follows that
$$
CL(\CP^n({\bf q});\Z)+1=SB(\CP^n({\bf q});\Z)=n+1.
$$
As a generalization of Fortune's theorem to  $(\CP^n({\bf
q}),\omega_{\rm FS}^{\bf q})$ we have:

\begin{theorem}\label{th:1.1}
If $h:\CP^n({\bf q})\times [0,1]\to\R$ is a $C^1$ time-dependent
Hamiltonian on $\CP^n({\bf q})$, then the time one map $\phi_1$ of
$X_{h_t}$ has at least $n+1$ fixed points. That is, the Arnold
conjecture $(AC_1)$ holds on $(\CP^n({\bf q}),\omega_{\rm FS}^{\bf
q})$.
\end{theorem}

No doubt the Oh's main result in \cite{Oh1} can be directly
generalized to $T^{2n}\times\CP^k({\bf q})$.

As a natural generalization $\RP^n\subset\CP^n$ we introduce a
suborbifold $\RP^n({\bf q})\subset\CP^n({\bf q})$ as follows:
$$
\RP^n({\bf q})=(\R^{n+1}\setminus\{0\})/\R^\ast,
$$
called {\bf real projective space of weight} ${\bf q}$. Here the
action of $\R^\ast$ on $\R^{n+1}\setminus\{0\}$ is still defined
by (\ref{e:1.1}). The isotropy group at any point ${\bf
x}\in\R^{n+1}\setminus\{0\}$ is given by $(\R^\ast)_{\bf
x}:=\cap_{x_j\ne 0}G_{q_j}\cap\R$, where $G_{q_j}=\{e^{2i\pi
k/q_j}\,|\, j=0,\cdots, q_j-1\}$ is the group of $q_j$-th roots of
unity. Clearly, $(\R^\ast)_{\bf x}$ is a subgroup of
$\Z_2=\{1,-1\}$, and thus $\RP^n({\bf q})$ is an orbifold of
dimension $n$. Clearly,  we have an orbifold isomorphism
\begin{equation}\label{e:1.4}
\RP^n({\bf q})\cong (\R^{n+1}\cap S^{2n+1}({\bf q}))/\Z_2,
\end{equation}
 where the action of $\Z_2$ is induced by
(\ref{e:1.1}). From this it easily follows that $\RP^n({\bf q})$
{\bf is a manifold if and only if all integers $q_1,\cdots,
q_{n+1}$ are odd.} In this case $\RP^n({\bf q})$  is diffeomorphic
to $\RP^n$. Hence $CL(\RP^n({\bf q});\Z_2)=n$ and $SB(\RP^n({\bf
q});\Z_2)=n+1$. If $\CP^n({\bf q})\equiv  S^{2n+1}({\bf q})/S^1$
is as above, and $[{\bf z}]_q$ denotes the $S^1$-orbit of ${\bf
z}\in S^{2n+1}({\bf q})$, then $\RP^n({\bf q})$ is isomorphic to
the Lagrangian suborbifold
 $$
L:=\{[{\bf z}]_q\,|\, \exists {\bf w}\in [{\bf z}]_q,\;{\bf
w}={\bf u}+ i{\bf v},\;{\bf v}=0\}
$$
in $(\CP^n({\bf q}),\omega_{\rm FS}^{\bf q})$.  As a
generalization of a result due to Chang-Jiang \cite{ChJi} and
Givental \cite{Gi} we have:

\begin{theorem}\label{th:1.2}
Let $q_1,\cdots, q_{n+1}$ be all odd. Then for any Hamiltonian map
$\phi_1:\CP^n({\bf q})\to\CP^n({\bf q})$ as in
Theorem~\ref{th:1.1}, it holds that $\sharp(\phi_1(\RP^n({\bf
q}))\cap\RP^n({\bf q}))\ge n+1$. Namely, in this case $(AC_2)$
holds for $(\CP^n({\bf q}), \RP^n({\bf q}))$.
\end{theorem}

Let $r({\bf q}):=\sharp\{q_i\,|\, q_i\in 2\Z+1\}$. Note that
$\R^{n+1}\cap S^{2n+1}({\bf q})$ is always diffeomorphic to $S^n$.
If $r({\bf q})=0$, by (\ref{e:1.4}) $\RP^n({\bf
q})\equiv\R^{n+1}\cap S^{2n+1}({\bf q})$ and hence  $SB(\RP^n({\bf
q}))=2$ and $CL(\RP^n({\bf q}))=1$.  If  $1\le r({\bf q})\le n$
then $\RP^n({\bf q})$ is an orbifold, not a manifold. In fact,
topologically $\RP^n({\bf q})$ is a $(n+1-r({\bf q}))$-fold
unreduced suspension of $\RP^{r({\bf q})-1}$, i.e., $\RP^n({\bf
q})=\Sigma^{(n+1-r({\bf q}))}(\RP^{r({\bf q})-1})$. It follows
that $SB(\RP^n({\bf q});\Z_2)=r({\bf q})$ and
$$
CL(\RP^n({\bf q});\Z_2)=1\;{\rm if}\;2\le r({\bf q})\le n,\;
 CL(\RP^n({\bf q});\Z_2)=0\;{\rm if}\; r({\bf q})=1.
$$
As an example we take ${\bf q}=(2,2,3)$, then $\RP^3((2,2,3))$ is
homeomorphic to the unit disk $D^2=\{x_1^2+ x_2^2\le 1\}$, and
thus for any group $G$ it holds that
$$
CL(\RP^3((2,2,3));G)=0\quad{\rm and}\quad SB(\RP^3((2,2,3));G)=1.
$$

 Theorem~\ref{th:1.1} and Theorem~\ref{th:1.2} were observed when
 author lectured a course of Variational methods for graduates
 from March to June 2005. Though their proofs can be obtained
by suitably changing ones in \cite{Fo} and \cite{ChJi}
respectively I am also to give main proof steps and necessary
changes.

\section{Proofs of Theorems}\setcounter{equation}{0}

\subsection{Proof of Theorem~\ref{th:1.1}}

 {\bf Later when talking the $S^1$-action on $\C^{n+1}$ or
$S^{2n+1}({\bf q})$ we always mean one given by (\ref{e:1.2})
without special statements}.  Denote by $B^{2n+2}({\bf q})=\{{\bf
z}\in\C^{n+1}\,|\, K_{\bf q}({\bf z})\le 1/2\}$. It is a compact
convex set containing the origin as an interior point and has
boundary $S^{2n+1}({\bf q})$. Since each nonzero ${\bf
z}\in\C^{n+1}$ can be uniquely expressed as ${\bf z}=r_{\bf z}{\bf
z}'$, $r_{\bf z}>0$ and ${\bf z}'\in S^{2n+1}({\bf q})$, for a
smooth family of functions $h_t:\CP^n({\bf q})\to\R$, $0\le t\le
1$, we can uniquely define a smooth family of functions
$H_t:\C^{n+1}\to\R$ by
\begin{equation}\label{e:2.1}
H_t({\bf z})=r_{\bf z}^2H_t({\bf z}')=r_{\bf z}^2h_t(\Pi({\bf
z}'))\quad{\rm and}\quad H_t(0)=0.
\end{equation}
Clearly, each $H_t$ is invariant under the action in
(\ref{e:1.2}), and positive homogeneous of degree two and
restricts to $h_t\circ\Pi$ on $S^{2n+1}({\bf q})$.  By the
standard symplectic reduction theory, cf., [McSa] and [LeTo], for
any constant $\lambda$ it is easily checked that
$\Pi_\ast(X_{H_t}+\lambda X_{K_{\bf q}})=X_{h_t}$. If $z:[0,
1]\to\CP^n({\bf q})$ satisfies $\dot z(t)=X_{h_t}(z(t))$ and
$z(0)=z(1)$ then there must exist a $\tilde z:[0, 1]\to
S^{2n+1}({\bf q})$ to satisfy $\dot{\tilde z}(t)=X_{H_t}(\tilde
z(t))$ and $z(t)=\Pi(\tilde z(t))$ for any $0\le t\le 1$, and
hence there is some $s\in [0, 2\pi)$ such that $\tilde z(0)=A_{s+
2k\pi}^{\bf q}\tilde z(1)$ for any $k\in\Z$. Define $\tilde
u_k(t)=A^{\bf q}_{t(s+ 2k\pi)}\tilde z(t)$ for $t\in [0,1]$. Then
$\tilde u_k(0)=\tilde u_k(1)$ and $\dot{\tilde u}_k(t)=X_{H_t+
\lambda K_{\bf q}}(\tilde u_k(t))$ for any $t\in [0,1]$ and
$\lambda\in s+ 2\pi\Z$; see the proof of Proposition 1 on page 144
in \cite{HoZe}. Conversely, if $\tilde z:[0, 1]\to\C^{n+1}$
satisfies $\tilde z(0)=\tilde z(1)\in S^{2n+1}({\bf q})$ and
$\dot{\tilde z}(t)=X_{H_t+ \lambda K_{\bf q}}(\tilde z(t))$ for
some $\lambda\in\R$ and any $t\in [0,1]$, then $\tilde
z([0,1])\subset S^{2n+1}({\bf q})$ and $z=\Pi\circ\tilde z:[0,
1]\to\CP^n({\bf q})$ satisfies $\dot z(t)=X_{h_t}(z(t))$ and
$z(0)=z(1)$; moreover for two such pairs $(\tilde z_1,\lambda_1)$
and $(\tilde z_2,\lambda_2)$, $z_1=\Pi\circ\tilde
z_1=\Pi\circ\tilde z_2=z_2$ implies $\lambda_1-\lambda_2\in
2\pi\Z$.
   Hence each
closed integral curve $z$ of $X_{h_t}$ on $\CP^n({\bf q})$
corresponds to a family
$$
\Omega_z:=\Bigl\{(u_k, s+
2k\pi)\,\Bigm|\!\!\!\!\!\left.\begin{array}{ll}
 &{\dot u}_k(t)=X_{H_t+ s+ 2k\pi K_{\bf
q}}(u_k(t)),\,u_k(0)=u_k(1) \\
&{\rm and}\;\Pi\circ u_k=\Pi\circ u_0 \,\forall k\in\Z
\end{array}\right.\!\!\!\!\Bigr\}.
$$
Clearly, the family $\Omega_z$ is diffeomorphic to $S^1\times
(2\pi\Z)$, and different families correspond to different fixed
points of $\phi_1$. So it suffice to prove:

\begin{equation}\label{e:2.2}
{\rm There}\;{\rm are}\;{\rm always}\;{\rm at}\;{\rm
least}\;(n+1)\;{\rm distinct}\;{\rm families}\;{\rm as}\;\Omega_z.
\end{equation}

In order to transfer it into a variational problem let
$Z:=L^2(\R/\Z,\C^{n+1})$ and
$$
X=\Bigl\{u=\sum_{k\in\Z}u_k\exp(2\pi ikt)\in Z\,\Bigm|\,
|u_0|^2+\sum_{k\in\Z}|k||u_k|^2<\infty\Bigr\}.
$$
Both carry respectively complete Hermitian inner products
$$
(u,v)=\sum_{k\in\Z}(u_k, v_k)_{\C^{n+1}}\quad{\rm and}\quad
\langle u,v\rangle=(u_0,v_0)_{\C^{n+1}}+
\sum_{k\in\Z}|k|(u_k,v_k)_{\C^{n+1}},
$$
where
$(\;,\;)_{\C^{n+1}}$ is the standard Hermitian inner-product on
$\C^{n+1}$. Let $ |u|=(u,u)^{\frac{1}{2}}$ and $\|u\|=\langle u,
u\rangle^{\frac{1}{2}}$ be the corresponding norms.
 Denote by $X^+=\{u\in X\,|\, u_k=0\forall k\le 0\}$,
$X^-=\{u\in X\,|\, u_k=0\forall k\ge 0\}$ and $X^0=\{u\in X\,|\,
u_k=0\forall k\ne0\}$. Then the natural splitting of $X$,
$X=X^+\oplus X^0\oplus X^-$, is an orthogonal decomposition of $X$
for the scalar products $\langle\;,\;\rangle$ and $(\;,\;)$. Let
$P_+$, $P_0$ and $P_-$ be the corresponding orthogonal
projections. Consider the densely defined self-adjoint linear
operator ${\cal L}:Z\supset D({\cal L})\to Z$ given by ${\cal
L}u=-i\dot u$. Then $\sigma({\cal L})=2\pi\Z$ and each $2\pi k$
has multiplicity $n+1$; moreover, ${\rm Ker}({\cal
L})\cong\C^{n+1}$ and normalized eigenvectors corresponding to
$2\pi k\in 2\pi\Z$ are $\phi_{k,j}=e^{2\pi ikt}\varepsilon_j$,
where $\varepsilon_1,\cdots,\varepsilon_{n+1}$ are the canonical
basis in $\C^{n+1}$.   Define  a compact self-adjoint linear
operator $L:X\to X$ by
\begin{equation}\label{e:2.3}
L(u)=2\pi (u^+-u^-)=2\pi\sum_{k\in\Z}ku_k,
\end{equation}
Clearly,  it is an extension of ${\cal L}$ to $X$ since $X$ can be
compactly embedded in $Z$. Consider a $C^1$ functional
$\Phi:X\to\R$ by
$$
\Phi(u)=\frac{1}{2}\langle Lu, u\rangle=\pi [\|u^+\|^2-\|u^-\|^2].
$$
For a time-dependent Hamiltonian $G_t$ on $\C^{n+1}$ we also
define ${\cal G}:X\to\R$ by
$$
{\cal G}(u)=\int^1_0G_t(u(t))dt.
$$
Set $\Phi_G=\Phi-{\cal G}$. As in Proposition 2.1 in [Fo] ones can
easily prove that $1$-periodic orbits of $X_{G_t}$ correspond to
critical points of $\Phi_G$ in a one-to-one way. The $S^1$-action
in (\ref{e:1.2}) induces an orthogonal $S^1$-representation
$\{T_s\}_{s\in S^1}$ on $X$ as follows:
\begin{equation}\label{e:2.4}
T_s(u)=\sum_{k\in\Z}(A_s^{\bf q}u_k)\exp(2\pi kt)
\end{equation}
if $u=\sum_{k\in\Z}u_k\exp(2\pi kt)$. ({\bf When saying
$S^1$-action on $X$ below we always mean this $S^1$-action without
special statements}).  The representation also preserve the
orthogonal splitting $X=X^+\oplus X^0\oplus X^-$. So if each $G_t$
is $S^1$-invariant then the functional ${\cal G}$ and thus
$\Phi_G$ is invariant with respect to the $S^1$-representation
$\{T_s\}_{s\in S^1}$.

Taking $H_t=G_t+ K_{\bf q}$ we can write
$$
\Phi_{H+\lambda K_{\bf q}}(u)=\Phi_H-\lambda {\cal K}_{\bf
q}(u),\quad{\rm where}\quad {\cal K}_{\bf q}(u)=\int^1_0K_{\bf
q}(u(t))dt.
$$
They are all $C^1$-smooth on $X$, and $S({\bf q}):={\cal K}_{\bf
q}^{-1}(1)$ is a Hilbert manifold. Since both $H_t$ and $K_{\bf
q}$ are positive homogeneous of degree two, by Lagrange multiplier
theorem the critical points of $\Phi_{H+\lambda K_{\bf q}}$ are in
a one-to-one correspondence with critical points of $\Phi_H$
constrained to $S({\bf q})$. Precisely speaking, if $u$ is a
critical point of $\Phi_H|_{S({\bf q})}$, then it is a critical
point of $\Phi_{H+\lambda K_{\bf q}}$ in $X$ with
$\lambda=\Phi_H(u)$; conversely, if $u$ is a critical point of
$\Phi_{H+\lambda K_{\bf q}}$ in $X$ and also sits in $S({\bf q})$,
 then $u$ is a critical point of $\Phi_H|_{S({\bf q})}$ with critical value
$\lambda$. Let $u_1, u_2$ be two critical points $\Phi_H|_{S({\bf
q})}$ with corresponding critical values $\lambda_1=\Phi_H(u_1)$
and $\lambda_2=\Phi_H(u_2)$. If $\lambda_1-\lambda_2\notin 2\pi\Z$
then $u_1$ and $u_2$ correspond to two geometrical different
$1$-periodic orbits of $X_{h_t}$ on $\CP^n({\bf q})$. Therefore
(\ref{e:2.2}) is reduced to prove:

\begin{theorem}\label{th:2.1}
There are always at least $(n+1)$ critical values of
$\Phi_H|_{S({\bf q})}$, $\lambda_i$, $i=1,\cdots, n+1$, such that
$\lambda_i-\lambda_j\notin 2\pi\Z$ for any $i\ne j$.
\end{theorem}

Clearly, we can always assume $h_t\ge 0$ and thus $H_t\ge 0$. By
(\ref{e:2.1}) ones have
$$
H_t({\bf z})\le 2\max_{(x,t)\in\CP^n({\bf q})\times I}h_t(x)
K_{\bf q}({\bf z})\le 2M K_{\bf q}({\bf z})
$$
for any ${\bf z}\in\C^{n+1}$ and $t\in\R$. Here
$M:=\sup_{(x,t)\in\CP^n({\bf q})\times I}h_t(x)$ and $I=[0,1]$. It
immediately gives the first claim in the following lemma.

\begin{lemma}\label{lem:2.2}
The functional ${\cal H}:X\to\R$ defined by ${\cal
H}(u)=\int^1_0H_t(u(t))dt$, satisfies:\\
 (i) $0\le{\cal
H}(u)\le M=2\max_{(x,t)\in\CP^n({\bf q})\times I}h_t(x)$ for any
$u\in{\cal K}_{\bf q}^{-1}(1)$;\\
(ii) $\nabla{\cal H}:X\to X$ is compact and equivariant, i.e.,
$\nabla{\cal H}(T_su)=T_s\nabla{\cal H}(u)$ for any $u\in X$ and
$s\in\R/2\pi\Z$.
\end{lemma}

The second properties is Proposition 2.3 in [Fo].

\begin{lemma}\label{lem:2.3}
The operator $X\to X, u\mapsto \nabla{\cal K}_{\bf q}(u)$ is
linear and bounded, and also respects the splitting of
$X=X^+\oplus X^0\oplus X^-$.
\end{lemma}

\noindent{\it Proof.}\quad It only need to check the final claim.
Let $u(t)=(u^{(1)}(t),\cdots, u^{(n+1)}(t))$ and
$u_k=(u^{(1)}_k,\cdots, u^{(n+1)}_k)$ for $k\in\Z$. Since
\begin{eqnarray*}
\langle\nabla{\cal K}_{\bf q}(u), v\rangle=d{\cal K}_{\bf
q}(u)(v)=&&\!\!\!\!\!\!\!\!\!\frac{1}{2}\int^1_0\sum^{n+1}_{j=1}q_j[(u^{(j)}(t),
v^{(j)}(t))+ (v^{(j)}(t), u^{(j)}(t))]dt\\
=&&\!\!\!\!\!\!\!\!\!\frac{1}{2}\sum^{n+1}_{j=1}\sum_{k\in\Z}q_j(u_k^{(j)}\bar
v_k^{(j)}+ v_k^{(j)}\bar u_k^{(j)})
\end{eqnarray*}
for any $u,v\in X$,  by the definitions of $X^+$, $X^0$ and $X^-$,
it easily follows that $\nabla{\cal K}_{\bf q}(X^+)\subset X^+$,
$\nabla{\cal K}_{\bf q}(X^0)\subset X^0$ and $\nabla{\cal K}_{\bf
q}(X^-)\subset X^-$. \hfill$\Box$\vspace{2mm}

Note that for any $u\in X$,
$$
\frac{1}{2}|u|^2\le\frac{\min_iq_i}{2}|u|^2\le{\cal K}_{\bf
q}(u)\le\frac{\max_iq_i}{2}|u|^2\le\frac{|{\bf q}|}{2}|u|^2.
$$
Carefully checking the proof of Proposition 3.2 in [Fo] ones can
easily use Lemma~\ref{lem:2.2} and Lemma~\ref{lem:2.3} to obtain:

\begin{lemma}\label{lem:2.4}
$\Phi_H|_{S({\bf q})}$ satisfies the Palais-Smale condition.
\end{lemma}

Note that $\Phi_H=\Phi-{\cal H}$. For $c\in\R$ and $\delta>0$ let
\begin{eqnarray*}
&&K_c:=\{u\in S({\bf q})\,|\, \Phi_H(u)=c,\; d(\Phi_H|_{S({\bf q})})(u)=0\},\\
&&N_\delta(K_c):=\{u\in S({\bf q})\,|\, d(u, K_c)=\inf_{v\in K_c}\|u-v\|<\delta\},\\
&&(\Phi_H|_{S({\bf q})})^c:=\{u\in S({\bf q})\,|\, \Phi_H(u)\le
c\}.
\end{eqnarray*}

Slightly changing the proof of Theorem 3.1 in [Fo] we can get:

\begin{lemma}\label{lem:2.5}
For any given $\delta>0$ and $c\in\R$ there exists an
$\varepsilon>0$ and an equivariant homeomorphism
$\eta:X\setminus\{0\}\to X\setminus\{0\}$ such that:
\begin{description}
\item[(i)] $\eta|_{S({\bf q})}$ is an equivariant homeomorphism to
$S({\bf q})$,

\item[(ii)] $\eta\bigl(S({\bf q})\cap(\Phi_H|_{S({\bf
q})})^{c+\varepsilon}\setminus N_\delta(K_c)\bigr)\subset S({\bf
q})\cap(\Phi_H|_{S({\bf q})})^{c-\varepsilon}$,
 \item[(iii)] $\eta(u)=Bu+
K(u)$, where $B:X\to X$ is an equivariant linear isomorphism of
the form $\exp(-tL)$ for some $t>0$ and $K:X\setminus\{0\}\to X$
is compact.
\end{description}
\end{lemma}

Note that the fixed point set of the $S^1$-action defined by
(\ref{e:2.4}),
$$
{\rm Fix}(\{T_s\}_{s\in S^1}):=\{u\in X\,|\, T_s(u)=u\forall
s\in\R/2\pi\Z\}=\{0\}.
$$
Let ${\cal A}$ be a family of all closed and $S^1$-invariant
subset $S\subset X\setminus\{0\}$, and
$$
\Lambda=\{h\in C^0(X, X)\,|\, h\circ T_s=T_s\circ h\,\forall
s\in\R/2\pi\Z\}.
$$
  Benci's index
[Be] is a map $\tau:{\cal A}\to\N\cup\{0,+\infty\}$ defined by
$$
\tau(S)=\left\{
\begin{array}{ll}
\min\{m\in\N\,|\, \exists \phi\in C^0(S,
\C^m\setminus\{0\}),\;\exists k\in\N:\\
\qquad \phi (T_su)=e^{2\pi
 iks}\phi (u)\,\forall (u, s)\in S\times\R/2\pi\Z\}, \qquad{\rm
 if}\;\{...\}\ne\emptyset,\\
 +\infty,\hspace{72mm}   {\rm if}\;\{...\}=\emptyset,
\end{array}
\right.
$$
and $\tau(\emptyset)=0$. For properties of the index $\tau$ see
Proposition 2.9 in \cite{Be}.

Let $\{R_s\}_{s\in\R/2\pi\Z}$ be an $S^1$-representation on $\C^k$
with $0$ as the only fixed point, and
$$
E^+=X^+,\quad E^-=X^-\oplus X^0,\quad P_{E^+}=P^+,\quad
P_{E^-}=P_-\oplus P_0.
$$
Then we get a $S^1$-representation $\{(T\oplus R)_s=T_s\oplus
R_s\}_{s\in\R/2\pi\Z}$ by
\begin{equation}\label{e:2.5}
(T\oplus R)_s(u\oplus x)=(T_su)\oplus(R_sx)
\end{equation}
for any $u\oplus x\in E^-\oplus\C^k$ and $s\in\R/2\pi\Z$, which
has $0$ as the only fixed point. In [BLMR] a relative index
(relative to $E^+$)
\begin{equation}\label{e:2.6}
\gamma(S)\in\N\cup\{0,+\infty\}\;{\rm of}\;{\rm a}\;{\rm
nonempty}\;{\rm set}\; S\in {\cal A}
\end{equation}
was defined as the minimum $m\in\N$ for which there is an
$S^1$-representation $\{R_s\}_{s\in\R/2\pi\Z}$ on $\C^m$ with the
fixed point set $\{0\}$, and an equivariant continuous map
$\phi:S\to (E^-\oplus\C^m)\setminus\{(0,0)\}$ with respect to
$S^1$-representations $\{R_s\}_{s\in\R/2\pi\Z}$ on $\C^m$ and
$\{(T\oplus R)_s=T_s\oplus R_s\}_{s\in\R/2\pi\Z}$ on
$E^-\oplus\C^m$ such that $P_{E^-}\circ\phi^-=P_{E^-}+ K$ for some
compact map $K:S\to E^-$(i.e, a continuous map which maps bounded
subsets in $S$ into relatively compact subsets in $E^-$). Here
$\phi^-$ is the $E^-$--component of $\phi$. As before, if no such
$m$ exist $\gamma(S)$ is defined as $+\infty$. Moreover,
$\gamma(\emptyset)=0$.

\begin{lemma}\label{lem:2.6}The relative index $\gamma$
satisfies:\\
{\rm (i)} $\gamma(S\cup R)\le\gamma(S)+ \gamma(R)$ for any $S,
R\in{\cal
A}$.\\
{\rm (ii)} If an equivariant isomorphism $h:X\to X$ leaves $E^-$
invariant, then for any $S\in{\cal A}$ it holds that $h(S)\in{\cal
A}$ and $\gamma(h(S))=\gamma(S)$.\\
{\rm (iii)} If $\gamma(S)\ge m$ and $E^+=F_1\oplus F_2$, where
$F_i$, $i=1,2$ are $S^1$-invariant and $\dim_{\C}F_1<m$,  then
$S\cap
F_2\ne\emptyset$.\\
{\rm (iv)} If $F\subset E^+$ is a complex $k$-dimensional
invariant subspace and $S(F,r)=\{u\in E^-\oplus F\,|\,{\cal
K}_{\bf q}(u)=r\}$, then
$\gamma(S(F,r))=k$.\\
{\rm (v)} For $S, R\in{\cal A}$, if there exists a continuous
bounded map $\varphi:S\to R$ such that
$P_{E^-}\circ\varphi=P_{E^-}+ K$
for some compact map $K:S\to E^-$, then $\gamma(S)\le\gamma(R)$.\\
{\rm (vi)} For $S, R\in{\cal A}$, if $\gamma(R)<\infty$ then
$\overline{S\setminus R}\in{\cal A}$ and
$\gamma(\overline{S\setminus R})\ge\gamma(S)-\tau(R)$.
\end{lemma}

\noindent{\it Proof.}  The proofs of these properties  can be
found in [BLMR] and [Be]. Ones only need to note that (iv) can be
proved by slightly changing the proof of Proposition 2.10 in
[BLMR]. In the present case we shall obtain a  map from the
elliptic sphere in $\C^n\times\C^l\times\C^k$ into
$\C^n\times\C^i\times\C^j$ with $j<k$ which is equivariant (with
respect to our $S^1$-action as the above $\{(T\oplus
R)_s\}_{s\in\R/2\pi\Z}$) and leaves $\C^l$, the fixed-point set,
invariant. It is not hard to check that the $S^1$-version of the
Borsuk-Ulam theorem due to [FHR] can be still used to get the
desired result.\hfill$\Box$\vspace{2mm}

For each $m\in\N$ let $\Gamma_m(S({\bf q}))=\{S\in{\cal A}\,|\,
S\subset S({\bf q}),\,\gamma(S)\ge m\}$ and
$$
c_m=\inf_{S\in\Gamma_m(S({\bf q}))}\sup_{u\in S}\Phi_H(u).
$$
For $j=1,\cdots, n+1$, let $\C^{n+1}_j\subset\C^{n+1}$ be the
complex $1$-dimensional subspace consisting of $w\in\C^{n+1}$
whose $k$-components are zero for $k\ne j$. Then for each
$m\in\N$,
$$
F_{m,j}=\oplus^m_{k=1}\C^{n+1}_j\exp(2\pi ikt)\subset E^+
$$
is a complex $m$-dimensional invariant subspace. By
Lemma~\ref{lem:2.6}(iv), $S(F_{m,j},r)\in\Gamma_m(S({\bf q}))$ and
$P_{E^+}(S(F_{m,j},r))$ is also compact. So it follows from
Lemma~\ref{lem:2.2}(i) that
\begin{eqnarray*}
 c_m\!\!\!\!\!\!&&\le\sup_{u\in S(F_{m,j},
1)}[\pi(\|u^+\|^2-\|u^-\|^2)-{\cal H}(u)]\\
&&\le\sup_{u\in S(F_{m,j},
1)}[\pi(\|u^+\|^2-\|u^-\|^2)]+ M\\
&&\le\sup_{u\in S(F_{m,j}, 1)}\pi\|u^+\|^2+ M\\
&&\le\sup_{u\in P_{E^+}(S(F_{m,j}, 1))}\pi\|u\|^2+ M<+\infty
\end{eqnarray*}
since $P_{E^+}(S(F_{m,j},1))$ is compact.
 Using Lemmas 2.4-2.6 the standard minimax arguments lead to:

\begin{theorem}\label{th:2.7}
$c_1\le c_2\le\cdots <+\infty$. If $c_m>-\infty$ then it is a
critical value of $\Phi_H|_{S({\bf q})}$. Moreover, if
$c_m=c_{m+1}=\cdots =c_{m+k}=c$ is finite then $\gamma(K_c)\ge k$.
\end{theorem}

Since $\Phi=\Phi_H+ {\cal H}$, if we set
$d_m=\inf_{S\in\Gamma_m(S({\bf q}))}\sup_{u\in S}\Phi(u)$, it
follows from the proof of Proposition 3.5 in [Fo] that
\begin{equation}\label{e:2.7}
d_m-M\le c_m\le d_m\quad\forall m\in\N.
\end{equation}
Here $M$ is defined as in Lemma~\ref{lem:2.2}(i). Clearly, the
$d_m$ have the same properties as the $c_m$ in
Theorem~\ref{th:2.7}. In particular, if $d_m$ is finite, it is a
critical value of $\Phi|_{S({\bf q})}$ and thus is an eigenvalue
of the linear eigenvalue problem
\begin{equation}\label{e:2.8}
Lu=\mu\nabla{\cal K}_{\bf q}(u)\quad{\rm on}\;S({\bf q}).
\end{equation}
Without loss of generality we now assume that
\begin{equation}\label{e:2.9}
q_1\ge q_2\ge\cdots\ge q_{n+1}\quad{\rm and}\quad q_1\ge 2.
\end{equation}

\begin{lemma}\label{lem:2.8}
Under the assumption (\ref{e:2.9}),  (\ref{e:2.8}) has eigenvalues
$$
\begin{array}{lcr}
\lambda^\pm_1=\pm\frac{2\pi}{q_1},\qquad
\lambda_2^\pm=\pm\frac{2\pi}{q_2},\qquad \cdots,
\lambda_{n+1}^\pm=\pm\frac{2\pi}{q_{n+1}},\\
\lambda_{1+ k(n+1)}^\pm=\pm\frac{2(k+1)\pi}{q_1},\; \lambda_{2+
k(n+1)}^\pm=\pm\frac{2(k+1)\pi}{q_2}, \cdots,
\lambda_{(k+1)(n+1)}^\pm=\frac{2(k+1)\pi}{q_{n+1}},\\
\lambda_0=0\;{\rm with}\;{\rm multiplity}\; n+1, k=1, 2,\cdots.
\end{array}
$$
Moreover, all $\phi_{k,j}=e^{2\pi ikt}\varepsilon_j$, $k\in\Z$ and
$j=1,\cdots,n+1$, are still the corresponding eigenvectors.
\end{lemma}

\noindent{\it Proof}. Assume that $Lu=\mu \nabla{\cal K}_{\bf
q}(u)$ for some $u\in S({\bf q})$ and $\mu\in\R$.  Let
$u(t)=\sum_{k\in\Z}u_k\exp(2\pi ikt)$, $u_k=(u_k^{(1)},\cdots,
u_k^{(n+1)})$, and $v(t)=\sum_{k\in\Z}u_k\exp(2\pi ikt)$,
$v_k=(v_k^{(1)},\cdots, v_k^{(n+1)})$. As in the proof of
Lemma~\ref{lem:2.3}, $\langle Lu, v\rangle=\mu \langle\nabla{\cal
K}_{\bf q}(u), v\rangle=\mu d{\cal K}_{\bf q}(u)(v)$ if and only
if
$$
2\pi\sum_{k\ne 0}\sum^{n+1}_{j=1}ku_k^{(j)}\bar v_k^{(j)}=
\frac{\mu}{2}\sum^{n+1}_{j=1}\sum_{k\in\Z}q_j(u_k^{(j)}\bar
v_k^{(j)}+ v_k^{(j)}\bar u_k^{(j)}).
$$
 The desired conclusions are easily derived from it.
\hfill$\Box$\vspace{2mm}

Denote by
\begin{equation}\label{e:2.10}
\cdots\le\mu_{-2}\le\mu_{-1}\le\mu_0=0<\mu_1\le\mu_2\le\cdots
\end{equation}
the sequence of eigenvalues of the eigenvalue problem
(\ref{e:2.8}),  each repeated according multiplicity. Let
$\hat\phi_k$ be the eigenfunction corresponding to $\mu_k$ for
$k\in\Z$. By Lemma~\ref{lem:2.3} we can normalize the $\hat\phi_k$
so that $\langle\nabla{\cal K}_{\bf
q}(\hat\phi_k),\hat\phi_l\rangle=\delta_{kl}$, $k,l\in\Z$. Note
that each $\hat\phi_k$ is the normalization of some $\phi_{i,j}$,
and that $k>0$ if and only if $i>0$. Moreover, it is clear that
$\{\hat\phi_k\,|\, k\in\Z\}$ form a complete orthogonal system in
$X$. So for each $u\in X$ it holds that $u=\sum_{k\in\Z}\langle
u,\hat\phi_k\rangle\hat\phi_k$. Especially,
$u=\sum_{k\in\Z}\langle u,\hat\phi_k\rangle\hat\phi_k\in S({\bf
q})$ if and only if
\begin{equation}\label{e:2.11}
1=\langle\nabla{\cal K}_{\bf q}(\sum_{k\in\Z}\langle
u,\hat\phi_k\rangle\hat\phi_k), \sum_{k\in\Z}\langle
u,\hat\phi_k\rangle\hat\phi_k\rangle=\sum_{k\in\Z}|\langle
u,\hat\phi_k\rangle|^2.
\end{equation}
Furthermore, assume that $u\in S({\bf q})\cap {\rm
span}\{\hat\phi_l\,|\, l\le k\}$; then we have
\begin{eqnarray*}
\Phi(u)=\frac{1}{2}\langle Lu,
u\rangle\!\!\!\!\!\!\!\!\!&&=\frac{1}{2}\Bigl\langle
L\Bigl(\sum_{l\le k}\langle u,\hat\phi_l\rangle\hat\phi_l\Bigr),
\sum_{l\le k}\langle u,\hat\phi_l\rangle\hat\phi_l\Bigr\rangle\nonumber\\
&&=\frac{1}{2}\sum_{l\le k}\sum_{j\le k}\langle
u,\hat\phi_l\rangle\langle u,\hat\phi_j\rangle\langle L\hat\phi_l,
\hat\phi_j\rangle\nonumber\\
&&=\frac{1}{2}\sum_{l\le k}\sum_{j\le k}\langle
u,\hat\phi_l\rangle\langle u,\hat\phi_j\rangle \mu_l\langle
\nabla{\cal K}_{\bf q}(\hat\phi_l), \hat\phi_j\rangle\\
&&=\frac{1}{2}\sum_{l\le k}\mu_l|\langle
u,\hat\phi_l\rangle|^2\le\mu_k.
\end{eqnarray*}
The final step is because of (\ref{e:2.11}). Hence we get
\begin{equation}\label{e:2.12}
d_k\le\mu_k\quad\forall k\ge 1.
\end{equation}
On the other hand, since $\phi_j\in E^+$ for any $j>0$, for any
$S\in\Gamma_k(S({\bf q}))$ with $k\ge 2$ it follows from
Lemma~\ref{lem:2.6}(iii) that the intersection $S\cap {\rm
span}\{\hat\phi_j| j\ge k\}$ must be nonempty since $E^+= {\rm
span}\{\hat\phi_j| 1\le j\le k-1\}+ {\rm span}\{\hat\phi_j| j\ge
k\}$ is an orthogonal  decomposition of invariant subspaces. Let
$v\in S\cap{\rm span}\{\hat\phi_j| j\ge k\}$. Then
$$
 \sup_{u\in S}\Phi(u)\ge
\Phi(v)=\frac{1}{2}\sum^\infty_{l=k}\mu_l|\langle
u,\hat\phi_l\rangle|^2\ge\mu_k.
$$
This and (\ref{e:2.12}) together yield:

\begin{proposition}\label{prop:2.9}{\rm ([Prop.3.6, BLMR])}
 $d_m=\mu_m$ for any $m>1$.
\end{proposition}

\begin{remark}\label{rem:2.10}{\rm
In Proposition 3.6 of [BLMR] it was also claimed that $d_1=\mu_1$.
However the arguments in the second step of proof therein seem not
to be complete for $k=1$. Precisely, for a set $B\in\Gamma_1(G_1)$
with $\gamma_r(pB)=1$, I do not know how  their Corollary 2.9 is
used to derive $B\cap{\rm span}\{\phi_i|i\ge 1\}\ne\emptyset$.
From the proof of their Proposition 2.8 it is impossible  to
improve their condition ``$\dim F_1<k$'' to ``$\dim F_1\le k$''.}
\end{remark}

Now (\ref{e:2.7}), Lemma~\ref{lem:2.8} and
Proposition~\ref{prop:2.9} together yield
\begin{equation}\label{e:2.13}
\frac{2\pi}{q_1}\le\mu_2\le\frac{2\pi}{q_2}\quad{\rm and}\quad
\mu_m-M\le c_m\le \mu_m\;\forall m\ge 2.
\end{equation}

\noindent{\bf Proof of Theorem~\ref{th:2.1}}.\quad Let $t_0\ge 1$
be the smallest integer such that $M\le 2t_0\pi$. By
(\ref{e:2.9}), $\max_iq_i=q_1\ge 2$. Then $4\pi/q_1\le 2\pi$ and
$$
2\pi/q_2\le\cdots\le 2\pi/q_{n+1}\le 2\pi.
$$
These imply that there are at least $(n+1)$'s $\mu_k$ (counting
multiplicity) in the interval $(2\pi/q_1, 2\pi]$. Using
Lemma~\ref{lem:2.8} it is easily seen that for each integer
$s>t_0+1$ the interval $(2(1+t_0)\pi, 2s\pi]$ contains at least
$(s-t_0-1)(n+1)$'s $\mu_k$. Let them be
$$
\mu_{l+1}\le\mu_{l+2}\le\cdots\le \mu_{l+ (s-t_0-1)(n+1)},\;l\ge
1.
$$
Then by (\ref{e:2.13}), corresponding with them we have
\begin{equation}\label{e:2.14}
c_{l+1}\le c_{l+2}\le\cdots\le c_{l+ (s-t_0-1)(n+1)},\;l\ge 1.
\end{equation}
By (\ref{e:2.13}) ones easily derive that these are in the
interval $I_s=(2\pi, 2s\pi]$, and thus that they are critical
values of $\Phi_H|_{S({\bf q})}$.

If there are $c_k, c_{k'}\in I_s$, $k, k'\ge l+1$, $k\ne k'$ such
that $c_k=c_{k'}$, the conclusion is obvious. So we can assume:
$$
c_k\ne c_{k'},\; k\ge l+1,\, k'\ge l+1,\;{\rm if}\; c_k\in I_s,
c_{k'}\in I_s.
$$
Then (\ref{e:2.14}) shows that
\begin{equation}\label{e:2.15}
\sharp(\{c_k\,|\, k\ge l+1\}\cap I_s)\ge (s-t_0-1)(n+1).
\end{equation}
Two elements $c$ and $c'$ in $\{c_k\,|\, k\ge l+1\}\cap I_s$ is
said to be equivalent if $c-c'$ is an integer multiple of $2\pi$.
Denote by $N_s$ the number of the equivalent classes. Without loss
of generality we can assume that $N_s$ is finite. Then
\begin{equation}\label{e:2.16}
\sharp(\{c_k\,|\, k\ge l+1\}\cap I_s)\le N_s(s-1).
\end{equation}
Take $s>1$ so large that $t_0(n+1)/(s-1)<1$. Then (\ref{e:2.15})
and (\ref{e:2.16}) give
$$
N_s\ge (n+1)-\frac{t_0(n+1)}{s-1}
$$
and thus $N_s\ge n+1$. The desired result is proved.
\hfill$\Box$\vspace{2mm}

\subsection{Proof of Theorem~\ref{th:1.2}}

 The original problem is reduced to estimate the
number of distinct solutions of the following boundary value
problem:
\begin{equation}\label{e:2.17}
\dot z(t)=X_{h_t}(z(t)),\; z(0)=p_0\in L,\; z(1)=p_1\in L.
  \end{equation}
Let $H_t$ be defined by (\ref{e:2.1}). As in \cite{ChJi}, modify
$H_t$ outside some open neighborhood of $B^{2n+2}({\bf q})$ so
that $H_t$ is $C^1$-bounded, and then consider a boundary value
problem:
\begin{equation}\label{e:2.18}
\left.\begin{array}{ll}
 &\dot z(t)=X_{H_t}(z(t))+ \lambda X_{K_{\bf q}}(z(t)),\\
 &z(j)\in \R^{n+1}\cap S^{2n+1}({\bf q}),\; j=0, 1.
 \end{array}\right\}
 \end{equation}
 Since $\Pi_\ast(X_{H_t}+\lambda
X_{K_{\bf q}})=X_{h_t}$, with the similar arguments to Lemma 2.1
and Lemma 2.2 in \cite{ChJi} we easily get:

\begin{lemma}\label{lem:2.11}
Each solution $z$ of (\ref{e:2.18}) sits in $S^{2n+1}({\bf q})$,
and $u(t)=\Pi(A^{\bf q}_{-\lambda t}z(t))$ solves (\ref{e:2.17}).
Moreover, if $(z^1, \lambda_1)$ and $(z^2,\lambda_2)$ are two
solutions of (\ref{e:2.18}), then
$$
\Pi(A^{\bf q}_{-\lambda_1 t}z^1(t))=\Pi(A^{\bf q}_{-\lambda_2
t}z^2(t))\,\forall t\Longrightarrow \lambda_1=\lambda_2 ({\rm mod}
\pi).
$$
\end{lemma}

Consider an operator $A: L^2([0,1],\C^{n+1})\supset D(A)\to
L^2([0,1],\C^{n+1})$ given by $Au=-i\dot u$, where
$D(A)=\bigl\{z\in H^1([0, 1],\C^{n+1})\,|\,
z(0),\,z(1)\in\R^{n+1}\bigr\}$. It is self-adjoint, and
$\sigma(A)=\pi\Z$. Moreover, each eigenvalue $k\pi$
 has multiplicity $n+1$, and corresponding eigenspace is spanned by
 $\varphi_{k,j}=e^{\pi ikt}\varepsilon_j$,
 $j=1,\cdots,n+1$. According to the spectral decomposition
 $$
 L^2([0,1],\C^{n+1})=\oplus_{k\in\Z}{\rm span}\{\varphi_{k,
 j}\,|\, j=1,\cdots, n+1\},
 $$
the operator $A$ can be decomposed into the positive, zero and
negative parts: $A=A^+ + A^0- A^-$. Let us  denote $D(|A|^{1/2})$
by
$$
X=\Bigl\{u=\sum_{k\in\Z}u_k\exp(\pi ikt)\in L^2([0,
1],\C^{n+1})\,\Bigm|\,
|u_0|^2+\sum_{k\in\Z}|k||u_k|^2<\infty\Bigr\}.
$$
It is a Hilbert space with inner product $(u,v)_X=\sum_{k\in\Z}(1+
|k|\pi)(u_k,v_k)_{\C^{n+1}}$ and the corresponding norm
$\|u\|_X=(u,u)_X^{1/2}$. Let ${\cal K}_{\bf q}:X\to\R$ be defined
by ${\cal K}_{\bf q}(u)=\int^1_0K_{\bf q}(u(t))dt$. Introduce the
manifold $S({\bf q})=\{u\in X\,|\, {\cal K}_{\bf q}(u)=1\}$ and
define a functional $J_H:X\to\R$ by
$$
J_H(u)=\frac{1}{2}(Bu,u)_X-\int^1_0H_t(u(t))dt,
$$
where $B:X\to X$ is defined by $B(u)=\pi\sum_{k\in\Z}ku_k$.
Slightly changing the proof of Lemma 2.4  in \cite{ChJi} one can
get:

\begin{lemma}\label{lem:2.12}
If $z_0\in S({\bf q})$ is a critical point of $J_H|_{S({\bf q})}$
and $\lambda_0$ is the corresponding Lagrange multiplier, then
$(z_0,\lambda_0)$ solves (\ref{e:2.18}) and $J_H(z_0)=\lambda_0$.
\end{lemma}

 There is an obvious $\Z_2$-action induced by (\ref{e:1.2}),
\begin{equation}\label{e:2.19}
g\cdot u=(g^{q_1}u_1,\cdots, g^{q_{n+1}}u_{n+1})\;\forall
g\in\Z_2=\{1,-1\},
\end{equation}
under which $J_H$ is invariant. Thus $J_H$ can be viewed as a
functional on the quotient $P({\bf q}):=S({\bf q})/\Z_2$. {\bf
Since all $q_1,\cdots, q_{n+1}$ are odd}, the action in
(\ref{e:2.19}) is free on $X\setminus\{0\}$, and {\bf hence
$P({\bf q})$ is a Hilbert manifold}. By (\ref{e:2.19}),
$-u=(-1)\cdot u$ for any $u\in X$, and thus $J_H(-z)=J_H(z)$ for
any $z\in X$. Note that for a given $z\in\C^{n+1}$,
$z\varphi_{k,j}\in X$ if and only if $z\in\R^{n+1}$. We set
\begin{eqnarray*}
 &&X_m=\oplus_{|k|\le m}\oplus^{n+1}_{j=1}(\R\varphi_{k,
 j}),\quad P({\bf q})_m=(S({\bf q})\cap
 X_m)/\Z_2,\\
&&X_m^+=\oplus^m_{k=1}\oplus^{n+1}_{j=1}(\R\varphi_{k,
 j}),\quad P({\bf q})^+_m=(S({\bf q})\cap
 X^+_m)/\Z_2,
\end{eqnarray*}
and $J_m=J_H|_{P_m({\bf q})}$. Then $Cl(\cup^\infty_{m=0}X_m)=X$.
By the proof of Lemma 3.1 in \cite{ChJi} ones can easily get: {\it
$J_H$ satisfies $(PS)^\ast$ with respect to the smooth Hilbert
filtration of finite dimension $P_1({\bf q})\subset P_2({\bf
q})\subset\cdots\subset P_m({\bf q})\subset\cdots$ of $P({\bf
q})$; that is, for any sequence $u^{(m)}\in P_m({\bf q})$,
$m=1,2,\cdots$, if $\{J_m(u^{(m)})\}$ is bounded and
$\lim_{m\to\infty}dJ_m(u^{(m)})=0$, then $\{u^{(m)}\}$ has a
convergent subsequence}.

The key is  that  $P_m({\bf q})$ is diffeomorphic to
$\RP^{(2m+1)(n+1)-1}$. Almost repeating the arguments in
\cite{ChJi} ones can get the desired result.
\hfill$\Box$\vspace{2mm}

\subsection{Open questions and concluding remarks} \setcounter{equation}{0}

(i) If $1\le r({\bf q})\le n$, the fixed point set of the action
in (\ref{e:2.19})  is given by
$$
{\rm Fix}_{\Z_2}=\{u=(u_1,\cdots, u_{n+1})\in X\,|\, u_i=0\;{\rm
if}\;q_i\notin 2\Z\}.
$$
Both ${\rm Fix}_{\Z_2}$ and $X\setminus{\rm Fix}_{\Z_2}$ are
infinite dimension subspaces. In this case the above methods fail.
Will $(AC_2)$ hold in the cases $0\le r({\bf q})\le n$?

\noindent{(ii)} Hofer's method in \cite{Ho} seem to be able to
prove the following result:

\noindent{\it Let $(M,\omega)$ be a closed symplectic orbifold and
$M^{sm}$ be its smooth locus. If $L\subset (M^{sm}, \omega)$ is a
compact Lagrange submanifold without boundary satisfying $\pi_2(M,
L)=0$, then for any Hamiltonian map $\phi:M\to M$ it holds that
$\sharp(L\cap\phi(L))\ge CL(L;\Z_2)+1$. Here $CL(L;\Z_2)$ denotes
the $\Z_2$-cuplength of $L$.}

Even if $\pi_2(M,L)\ne 0$, but $L$ is monotone and its minimal
Maslov number $N_L\ge 2$ it is also possible to generalize some
results in \cite{Oh2} to the case that $(M,\omega)$ is a closed
symplectic orbifold and $L\subset M^{sm}$.

\end{document}